# Strong Valid Inequalities Identification for Mixed Integer Programming Problems

Asghar Moeini[1], Kate Smith-Miles

*School of Mathematics & Statistics, The University of Melbourne, Parkville 3010, Australia*

## Abstract

The characterization of strong valid inequalities for integer and mixed-integer programs is more of an artistic task than a systematic methodology, requiring inspiration that can sometimes be elusive. Frequently, this task is facilitated by somehow exploiting the structure of problems for devising strong valid inequalities. Subsequently, various mathematical techniques are utilized for proving that those inequalities, which are often easily shown to be valid, are indeed strong in the sense that they represent facets or other high dimensional faces. This paper develops a method to assist modelers in the challenge to devise strong valid inequalities. In each iteration, the proposed algorithm generates a valid inequality by solving a suitably constructed linear mixed integer program and applies some quality criteria in order to determine if it is a new strong valid inequality. To illustrate the proposed algorithm, a new Traveling Salesman Problem (TSP) formulation is developed based on a set of constraints already constructed in the context of the Hamiltonian Cycle Problem (HCP), and then the proposed algorithm is employed to derive a set of strong inequalities to tighten this TSP formulation. Finally, a comparison study between the relaxation of the new TSP formulation and that of a state-of-the-art TSP formulation is conducted. The computational study confirms the effectiveness of the devised inequalities due to the better quality of the relaxation provided by the new formulation.

**Keywords**: Integer programming, Linear programming, Combinatorial optimization, Facet-defining inequality, Valid inequality, Convex hull problem, Traveling salesman problem

## 1. Introduction

A natural approach to study a broad range of optimization problems arising in operations research is to express them as Integer Programming (IP) formulations of the form

$$\begin{aligned} \text{Min} \quad & c^T x \\ \text{s.t.} \quad & Ax = b \\ & Bx \geq d \\ & x_1, \dots, x_q \in \mathbb{Z} \\ & x_{q+1}, \dots, x_n \in \mathbb{R} \end{aligned} \qquad \text{(IP)}$$

where $x$ is a vector representing the decision variables, $c \in \mathbb{R}^n$ is a vector representing objective function coefficients, $A \in \mathbb{R}^{l \times n}$ and $B \in \mathbb{R}^{m \times n}$ are matrices representing the coefficients for equality and inequality constraints respectively, and $b \in \mathbb{R}^l$ and $d \in \mathbb{R}^m$ are vectors representing the right hand side values for those constraints. IP problems are known to be NP-hard [18], because of their integral variables. Therefore, it is logical to start their analysis by first investigating their LP-relaxation, namely the problems

$$\begin{array}{rl} \text{Min} & c^T x \\ \text{s.t.} & Ax = b \\ & Bx \geq d \\ & x \in \mathbb{R}^n \end{array} \qquad \text{(LP)}$$

[1] Corresponding Author Email: asghar.moeini@unimelb.edu.au

A critical question that arises is whether we can further reduce $\mathcal{X}$, the feasible space of LP, without removing any solutions lying in $S$, the feasible space of the IP problem. To tackle such a reduction task, we need to study the properties of polytopes associated with these problems.

It is well known that a polytope can be specified in two ways [33]. The V-representation defines a polytope with the set of its extreme points (vertices), while the H-representation defines a polytope through a set of valid inequalities. An important question in the integer programming context is how we can compute an H-representation of a polytope from its V-representation. That is, if $S$ is a finite set, then the convex hull of $S$, denoted by **conv**($S$), is a well-defined polytope and according to the results in [30], there is at least a finite set of valid inequalities that characterizes **conv**($S$). Such a linear description is called the H-representation of the set $S$ (see also [31]). Computing the H-representation for a given set of vertices $S$ is called the facet enumeration problem and has been studied by numerous researchers, see for example [3], [4], and [14]. There are also some software packages such as Polymake [15] and Porta [6] for polytope specification purposes. The performance and implementation aspects of the state-of-the-art facet enumeration algorithms is summarized in [1]. Note that in this manuscript, we utilize the term "H-representation" for the polyhedral description of any given instance of the IP problem, and we use the term "convex hull formulation" for the polyhedral description in the general form of the underlying problem.

Characterizing the convex hull formulation of an NP-hard problem is always heuristic in nature since if there existed such a description, which was also computationally tractable (i.e., there was a polynomial separation algorithm to recognize a violated inequality for any given point outside the polytope), then the original NP-hard problem could be solved in polynomial time [17]. Therefore, one may ask why we are still interested in finding valid inequalities to refine the convex hull formulations if they cannot be fully characterized? This is because in practice, even an insufficient description with a few computationally tractable families of inequalities is adequate to solve many problems. For example, see [13], in which branch-and-cut algorithms are utilized to solve 35 large-sized real-world instances of Asymmetric Traveling Salesman Problem.

The above discussion reveals the importance of identifying strong valid inequalities to approximate convex hull formulations. To the authors' knowledge, there is no formal method of doing so, and it is often achieved by using knowledge of the structure of the underlying problem. In particular, this means trying to guess a family of valid inequalities and then proving that they hold for all instances of the underlying problem. However, sometimes even a complicated structural analysis of a problem may not result in a useful constraint (equality/inequality). For example, in [12], an analysis of properties of Hamiltonian cycle matrices led to a set of constraints which were acknowledged as redundant constraints.

One popular approach for guessing valid inequalities is to obtain the H-representation of very small instances, and then try to recognize generalizable patterns. As an example of this approach we refer the reader to [5] in which a few new facets of a TSP formulation derived. This guessing process to devise strong valid inequalities often involves dealing with several issues. The most important issue with methods that are based on the H-representation is their dimensional limitations [2], because even H-representation of small instances often include an extremely large number of inequalities. For example, constructing the H-representation of 10-node TSP polytope requires about fifty billion facet-defining inequalities [21].

This limitation does not allow consideration of large enough instances and, therefore, the generated inequalities (H-representations) might not be generalizable, especially due to the fact that the structures of small polytope instances could be different from the structures of large polytope instances. Thus, the small H-representations are often inadequate to develop an intuition as to how the behaviours would extend for the general case (convex hull formulation). This implies that in order to study the polytope of an IP problem, we need to consider the H-representation of sufficiently large instances, and enough of them, to increase the chance of detecting generalizable inequalities.

Another issue with the guessing process is that the H-representation of a polytope can be variant in the sense that there could be several sets of constraints (equalities/inequalities) to represent the same polytope. This means the H-representation methods can sometimes result in a set of constraints in which the underlying patterns are hidden and beyond easy recognition, making generalization difficult.

Motivated by these issues, the present paper proposes a potential method to identify strong valid inequalities for IP models. The proposed method can also be employed in the context of extended formulations (e.g., see [7]) to strengthen an IP model by involving new decision variables and then devising new inequalities according to the relationships between the original and new decision variables.

The remainder of the paper is arranged as follows. Definitions and preliminaries are introduced in Section 2. The proposed method for generating strong valid inequalities is presented in Section 3. To illustrate the advantages of the proposed method, a new formulation for the TSP is developed in Section 4 and then in Section 5, the proposed method is employed to devise a set of strong inequalities for the new TSP formulation. Furthermore, a numerical comparison is given between the LP-relaxations of our tightened TSP formulation and one of the state-of-the-art TSP formulations to show its competitiveness. Given that the new TSP formulation was initially uncompetitive due to its loose LP-relaxation, the ability of the proposed method to tighten the LP-relaxation to the point of being competitive with the LP-relaxation of a state-of-the-art model shows the promising direction established by the proposed approach. Finally, conclusions are stated in Section 6.

## 2. Definitions and Preliminaries

This section sets up required definitions and preliminaries and offers a brief review of our previously published *Equality Constraint Augmenting* (ECA-method) [25] which plays a prominent role in establishing the proposed method for strong valid *inequality* identification.

A set $\mathcal{g} \subseteq \mathbb{R}^n$ is affine if the line through any two points in $\mathcal{g}$ is entirely contained in $\mathcal{g}$ (i.e., $\alpha s_1 + (1-\alpha)s_2 \in \mathcal{g}$ when $x_1, x_2 \in \mathcal{g}$ and $\alpha \in \mathbb{R}$). Generalizing this idea, we refer to a point of the form $\alpha_1 s_1 + \cdots + \alpha_N s_N$, where $\alpha_1 + \cdots + \alpha_N = 1$, as an affine combination of the points $s_1, \ldots s_N$ . It can easily be seen that an affine set contains every affine combination of its points. The set of all affine combinations of points in some set $S \in \mathbb{R}^n$ is called the affine hull of $S$, and is denoted by $\mathbf{aff}(S)$:

$$\mathbf{aff}(S) = \{\alpha_1 s_1 + \cdots + \alpha_N s_N | s_1, \ldots, s_N \in S, \alpha_1, \ldots, \alpha_N \in \mathbb{R}, \alpha_1 + \cdots + \alpha_N = 1\}.$$

The points $s_1, \ldots, s_N \in \mathbb{R}^n$ are affinely independent if the unique solution to $\alpha_1 s_1 + \cdots + \alpha_N s_N = 0$ is $\alpha_1, \ldots, \alpha_N = 0$. The dimension of an affine set $\mathcal{g}$, denoted by **dim**($\mathcal{g}$), is the maximum number of affinely independent points in $\mathcal{g}$ minus one. Analogously, the dimension of $\mathbf{conv}(S)$, the polytope defined as the convex hull of points $s_1, \ldots, s_N$, is equal to the dimension of the affine hull of these points.

It can easily be shown that (e.g. see [31]) if a polytope $P$ is constructed as the convex hull of points $s_1, \ldots, s_N$ then its dimension can be obtained as

$$\mathbf{dim}(P) = \mathbf{rank}([v_1 | v_2 | \ldots | v_{N-1}]),$$

where $v_{i-1} := s_i - s_1,\ i = 2, \ldots, N$. Furthermore, if the polytope $P \subseteq \mathbb{R}^n$ is given as a set of constraints, namely $P = \{Ax = b, Bx \geq d\}$ then it follows that

$$\mathbf{dim}(P) = n - \mathbf{rank}(A),$$

providing that $P$ is nonempty and the inequalities $Bx \geq d$ do not imply any additional equality constraints. We need the following definitions to formally continue our exposition.

**Definition 1.** The inequality $\pi^T x \leq \pi_0$ is a valid inequality for polytope $P$ if it is satisfied by all points in $P$, i.e., if and only if $\max\{\pi^T x \mid x \in P\} \leq \pi_0$.

**Definition 2.** If $\pi^T x \leq \pi_0$ is a valid inequality for a polytope $P$, then $f = \{x \in P | \pi^T x = \pi_0\}$, $f$ is called a face of $P$, and we say that $(\pi, \pi_0)$ represents $f$.

**Definition 3.** A face $F$ of polytope $P$ is called a facet of $P$ if **dim**($F$) = **dim**($P$) – 1.

**Definition 4.** A polytope $P$ is called full dimensional if there is no equality constraint which can be satisfied by all points $x \in P$.

We also use the following known results from linear algebra (e.g. see [32]).

**Theorem 1.** *If $F$ is a facet of $P$, there is some inequality that represents $F$.*

**Theorem 2.** *If polytope $P \subseteq \mathbb{R}^n$ and* **dim**($P$) = $n - d$*, where $d > 0$, then the maximal number of independent equality constraints that can be satisfied by points in P is $d$.*

**Theorem 3.** *If $f$ is a $(n-d)$-dimensional face of polytope $P \subseteq \mathbb{R}^n$, where $d > 1$, then $f$ lies in intersections of d facets of $P$.*

**Theorem 4.** *If $P$ is a full dimensional polytope, it has a unique minimal description $P = \{x \in \mathbb{R}^n \mid Bx \leq d\}$ where each inequality is unique to within a positive multiple.*

Suppose $S$ denotes the solution set of (IP). According to Theorem 2, if the dimension of **conv**($S$) is $d$ units less than the dimension of $\mathcal{X}$, this implies that there exist $d$ unidentified equality constraints. These new constraints can be extracted using the ECA-method [25]. The following section provides a review of this method.

## 2.1. Review of Equality Constraint Augmenting Method

Throughout this paper, we assume that the set of solution of (IP) is available. However, in practice we might need to use Lattice point enumeration methods to obtain such a set. We refer to [1] in which the performance of the state-of-the-art Lattice point enumeration algorithms is reported.

Let $S = \{s_i \in \mathbb{R}^n,\ i = 1,\dots,N\}$ denote the set of solutions of (IP) and define: $v_{i-1} := s_i - s_1,\ i = 2,\dots,N,$ and $V := [v_1|v_2|\dots|v_{N-1}]$ the $n \times N$ matrix whose columns are $v_i$'s. Ignoring the existing equality constraints, $A_0x = b_0$, the ECA-method finds a set of all equality constraints, $\tilde{A}x = \tilde{b}$ and then removes the ones captured by the already known equality constraints. Thus, only the equalities that refine the feasible space further are retained.

To extract all equality constraints, one needs to find a linear equation system, $\tilde{A}x = \tilde{b}$, which represents the affine hull of $S$ (**aff**($S$)), that is, $\{\tilde{A}x = \tilde{b} \mid x \in \mathbf{aff}(S)\}$. It is shown in [24] and [25] that a matrix $\tilde{A}$ and vector $\tilde{b}$ can be obtained as follows:

$$\tilde{A} := \mathbf{null}(V^T)^T, \qquad \tilde{b} := \tilde{A}s_1, \tag{*}$$

where **null**($K$) denotes a matrix, whose columns are a basis of the null space of matrix $K$, and $s_1$ is one of the solutions in the set $S$. The generated set of equalities do not depend on the choice of $s_1$ , thus $s_1$ can be chosen as any point in $S$. Note that the matrix $V^T$ is often large and sparse which allows exploitation of some efficient methods to obtain the basis of the null space [16]. Also, Matlab software is available [20] to find the basis of the null space of such matrices.

The generated equalities $\tilde{A}x = \tilde{b}$ provide a set of all (existing and unidentified) equality constraints. Therefore, these generated equalities must be consecutively checked to see if they are unidentified. That is, the first generated equality is appended to the existing equality constraints to see whether it increases the rank of the coefficient matrix by one. If so, the constraint is new (non-redundant), and the existing equalities will be updated by including this new equality. Repeating this process guarantees all newly added equalities are new. Note that, the order in which the new constraints are added effects which ones are retained. Thus, we sort the generated equalities ($\tilde{A}x = \tilde{b}$) based on the number of their non-zero coefficients, increasing the chance that the constraints to retain are easier for generalization.

Defining $K(i.)$ and $K(.i)$ as the $i^{\text{th}}$ row and $i^{\text{th}}$ column of matrix $K$ respectively, the ECA-method can be represented as Algorithm 1.

**Algorithm 1 (ECA-method):**

**Step 0.** Receive the inputs including existing equality constraints, $A_0x = b_0$ and set of solutions, $S = \{s_i \in \mathbb{R}^n,\ i = 1,\ldots,N\}$. Construct $V$ based on $S$. Set $i = 0$, $E = A_0$, and $j = 1$.

**Step 1.** Obtain $m_0 = \mathbf{Rank}(A_0)$, and $m = \mathbf{Rank}(V)$. Set $d = m - m_0$.

**Step 2.** Stop if $d = 0$ and return "There are no unidentified equality constraints", otherwise, return "There are $d$ unidentified equality constraints".

**Step 3.** Find all equality constraints, $\tilde{A}x = \tilde{b}$ based on the equations in (*) and rearrange the rows of the matrix $\tilde{A}$ based on sorting the number of their non-zero coefficients to construct the matrix $A_1$.

**Step 4.** Construct matrix $E = \left[\frac{E}{A_1(j.)}\right]$ [2]. If $\mathbf{Rank}(E) = m_0 + 1$, then update $i = i + 1$, and $Q(i.) = A_1(j.)$. Otherwise, update $E$ by eliminating its last row.

**Step 5.** Stop if $i = d$, and return the set of new equality constraints, $Qx = b_Q$, where $b_Q = Qs_1$. Otherwise, update $j = j + 1$ and go back to Step 4.

## 3. The proposed method to generate strong valid inequalities

In this section, we propose a potential approach to identify new strong valid inequalities for IP problems that can be included already with a partial set of (equality/inequality) constraints. We call this method *Inequality Constraint Augmenting method*, or ICA-method, for short. Note that this method can be employed to extract strong valid inequalities including faces of lower dimension than that of facets. However, we explain it here only for finding facet-defining inequalities.

The main idea of the ICA-method is intuitively derived from the well-known fact that faces of polytopes are also polytopes themselves. According to Definition 3, the dimension of facets of a polytope is one less than the dimension of that polytope. That is, if the maximal number of independent equalities that the points of a given polytope can satisfy is $m$, then the points on a facet of that polytope can satisfy exactly $m + 1$ independent equalities. Let $\pi^T x = \pi_0$ be the equality that points lying in a facet can satisfy but other points of the polytope cannot. Suppose $s_0$ is a point of the polytope that does not lie in the facet. Therefore, according to Theorem 1, if $\pi^T s_0 < \pi_0$ then $\pi^T x \leq \pi_0$ is the facet-defining inequality, and otherwise if $\pi^T x > \pi_0$, then $-\pi^T x \leq -\pi_0$ represents the facet-defining inequality.

If the set of solutions lying on a facet, $S_F$ is known, we can exploit the aforementioned ECA-method to find $(\pi, \pi_0)$. More precisely, we first generate all $m$ equalities (if they are not already identified) that the solutions in $S$ satisfy, called Set 1 and then generate all $m + 1$ equalities that the solutions in $S_F$ satisfy, called Set 2. Subsequently, we can check which equality in Set 2 cannot be written as a linear combination of those $m$ equalities in Set 1. In order to achieve this, equalities in Set 2 can be added one at a time to the $m$ equalities of Set 1 to check if the rank of the newly constructed coefficient matrix is increased by one. Note that only one equality among $m+1$ equalities of Set 2 increases the rank and therefore that equality provides $(\pi, \pi_0)$.

[2] $\left[\frac{E}{A_1(j.)}\right]$: Embedding $A_1(j.)$ to the end of $E$.

A natural question that arises is how we can find a set of solutions (vertices) that lie in a facet. This is a very important question since if we can find such sets of solutions then the ECA-method can find the associated facets. From basic linear algebra, if a polytope with $N$ extreme points (vertices) has dimension $k$, then at least $k$ vertices lie in each facet of that polytope. Hence, for an exhaustive search, we need to take $k$ vertices out of $N$ and check if they provide a facet. This implies that all $\binom{N}{k}$ possible ways need to be checked for extracting all facets, and therefore it is impractical.

Note that facets are valid inequalities which are often satisfied by numerous solutions. Therefore, to find a facet of a polytope constructed as the convex hull of vertices $s_1, \dots, s_N$, we can design a mixed integer program (§) to find a valid inequality, $\pi^T x \leq \pi_0$ with the maximal number of solutions satisfying it at equality.

$$\begin{array}{llll} \text{Max} & z = \sum_{i=1}^{N} \theta_i & & \\ \text{s.t.} & -\pi^T s_i + \pi_0 \leq (1-\theta_i)M & i = 1, \dots N & (1) \\ & \pi^T s_i - \pi_0 \leq (\theta_i - 1)\varepsilon & i = 1, \dots N & (2) \\ & \sum_{i=1}^{N} \theta_i \leq N-1 & & (3) \\ & \sum_{i=1}^{N} \theta_i \geq n-m & & (4) \\ & \theta \in \{0,1\} & i = 1, \dots N & (5) \end{array} \qquad (§)$$

where $S = \{s_i \in \mathbb{R}^n,\ i = 1, \dots, N\}$ is the set of solutions of (IP). $\theta_i$'s are binary decision variables associated with solutions such that $\theta_i$ is one if $s_i$ satisfies the valid inequality at equality ($\pi^T s_i = \pi_0$) and otherwise is zero. $m$ is the number of all independent equality constraints ($m = \mathbf{Rank}(\tilde{A})$ ), and $\varepsilon$ and $M$ are a small and large positive value respectively.

In (§), the objective function counts the number of solutions that satisfy the valid inequality at equality. Constraints (1) and (2) together ensure that if the solution $s_i$ lies in the inequality (i.e., $\pi^T s_i = \pi_0$) then $\theta_i = 1$ and otherwise if $\pi^T s_i \leq \pi_0 - \varepsilon$ then $\theta_i = 0$. Constraint (3) prevents that all solutions can satisfy the valid inequality at the equality form, because it then will be an equality constraint rather than a valid inequality. Constraint (4) is written as each facet of a $(n-m)$-dimensional polytope can be satisfied by at least $n-m$ solutions.

**Proposition 1**. *Solving (§) provides a facet of the polytope* **conv**($S$).

**Proof**. By contradiction suppose solving (§) gives a face, $f$ which has a lower dimension than that of facets. According to Theorem 3, $f$ is a subset of some facets and each of those facets not only are satisfied at equality by all the solutions that lie in $f$ but they also can be satisfied by some other solutions at equality. Therefore, the number of solutions that lie in $f$ are less than the number of solutions lying on those facets and this is obviously a contradiction with our assumption ($f$ is the optimal solution of (§)). This implies that the solution of (§) cannot be anything but a facet and it completes the proof.

■

Note that the facet found by solving (§) might not be a new valid inequality and could be implied by existing constraints. To check whether it is a new or redundant inequality, one can solve the following linear program.

$$\begin{aligned} \text{Max} \quad & z = \pi^T x \\ \text{s.t.} \quad & \tilde{A}x = \tilde{b} \\ & B_0 x \geq d_0 \end{aligned} \qquad (\dagger)$$

where $\tilde{A}x = \tilde{b}$ and $B_0 x \geq d_0$ represent the set of all equality constraints and existing inequality constraints respectively. Suppose the optimal objective value of (†) is $z^*$ and $z^* \leq \pi_0$, then valid inequality $\pi^T x \leq \pi_0$ is redundant.

Regardless of whether solving (§) leads to a new or redundant inequality, one might be interested to re-solve (§) to identify another valid inequality. So, how we can prevent revisiting the former generated valid inequality when re-solving (§)? Before answering this question in Proposition 2, let define $\theta^* = (\theta_1^*, \ldots, \theta_N^*)^T$ to be values in an optimal solution of (§).

**Proposition 2**. *Updating the set of constraints of* (§) *by adding the constraint* (‡)*, eliminates the former optimal solution (found valid inequality) from the feasible space.*

$$\sum_{i \in \{j | \theta_j^* = 1\}} \theta_i - \sum_{i \in \{j | \theta_j^* = 0\}} \theta_i \leq -1 + \sum_{i=1}^{N} \theta_i^* \qquad (\ddagger)$$

**Proof**. Constraint (‡) ensures that either at least one of the $\theta_i$'s variables with value one in $\theta^*$ must switch to zero or at least one of $\theta_i$'s with value zero in $\theta^*$ must switch to one and therefore the model cannot deliver the former optimal solution as a new optimal solution.

■

**Remark 1.** Note that solving (§) gives a facet-defining inequality which is often in a form that is hard to interpret and generalize for higher dimensional instances of the underlying problem. Moreover, we already studied the process of using algorithm 1 for a given set of solutions that lie in a facet to find the associated facet-defining inequality. We can utilize that process to reobtain the complicated appearance facet-defining inequalities in a simpler form. More precisely, we can consider the set of solutions lying in this facet, $S_F = \{s_i | \theta_i = 1, i = 1, \ldots N\}$ as the input for the that process to regenerate an inequality which is likely to be simpler for generalization. We shall demonstrate this approach in example 1.

Note that constraint (‡) can also be applied to avoid regenerating any of the existing inequalities ($A_0 x = b_0$) as follows.

**Remark 2.** Suppose $\pi^T x \leq \pi_0$ is an existing valid inequality. We then construct its associated $\theta^*$ vector, a vector whose $i^{\text{th}}$ entry has value one if $\pi^T s_i = \pi_0$ and otherwise has value zero. Subsequently, we write the constraint (‡) for $\theta^*$ and add it to (§) which ensures we do not regenerate $\pi^T x \leq \pi_0$ as a valid inequality. We shall refer to this constraint as deduplication constraint for $\pi^T x \leq \pi_0$. If the deduplication constraints are not added to (§), we risk revisiting existing facets in further iterations.

The ICA-method is now summarized in Algorithm 2 as follows.

**Algorithm 2 (ICA-method):**

**Step 0.** Receive the inputs including existing equality constraints, $A_0x = b_0$ and set of solutions, $S = \{s_i \in \mathbb{R}^n,\ i = 1,\dots,N\}$. Generate a set of all equality constraints, $\tilde{A}x = \tilde{b}$ (if they are not already identified) by using Algorithm 1.

**Step 1.** Construct the model (§). Update (§) by embedding the deduplication constraints written for existing inequalities based on the instruction in Remark 2.

**Step 2.** Solve (§) strop if it is infeasible. Otherwise, let $\theta_i^*$ denote the $\theta_i$ value in an optimal solution, and find $S_F = \{s_i \mid \theta_i^* = 1, i = 1, \dots N\}$.

**Step 3**. If $S_F$ is a facet (i.e., $\mathbf{dim}(\mathbf{conv}(S_F)) = n - \mathbf{rank}(\tilde{A}) - 1$), then go to the next step, otherwise go to Step 6.

**Step 4**. Utilize (†) to check whether it is a new facet. If 'yes' go to the next step, otherwise go to Step 6.

**Step 5**. Regenerate a simpler form of this inequality by the instruction explained in Remark 1.

**Step 6.** Update the model (§) by adding the deduplication constraint written for $\theta^*$. Also, update Constraint (3) of (§) to $\sum_{i=1}^{N} \theta_i \leq z^*$, where $z^*$ is the optimal value of the objective function of (§) and go back to step 2.

Note that the optimal value of the objective function of (§) will be either equal to or less than the optimal objective of the former iteration. Hence, we can update the right-hand side of Constraint (3) of (§) as mentioned in step 6 of the above algorithm.

**Remark 3.** Note that there are polytopes with facets such that a small number of extreme points can satisfy those facets at equality, while the same polytopes have faces of lower dimensions with larger numbers of vertices lying on these faces. Therefore, in some iterations of Algorithm 2, we might reach a face rather than a facet, and we can retain the faces with sufficiently large dimensions as strong valid inequalities. It follows that Algorithm 2, which is designed for identification of facet-defining inequalities, can be modified to extract strong valid inequalities. For doing so, one needs to define a threshold dimension for considering the inequalities associated with faces of higher dimension than that threshold as strong valid inequalities.

We now present a simple example to illustrate the implementation of the proposed ICA-method.

**Example 1.** Suppose we are given the three-dimensional Stasheff polytope in $\mathbb{R}^4$ whose extreme points are labelled as shown in Figure 1, and whose 4-dimensional coordinates for the labelled extreme points of the polytope are defined in Table 1.

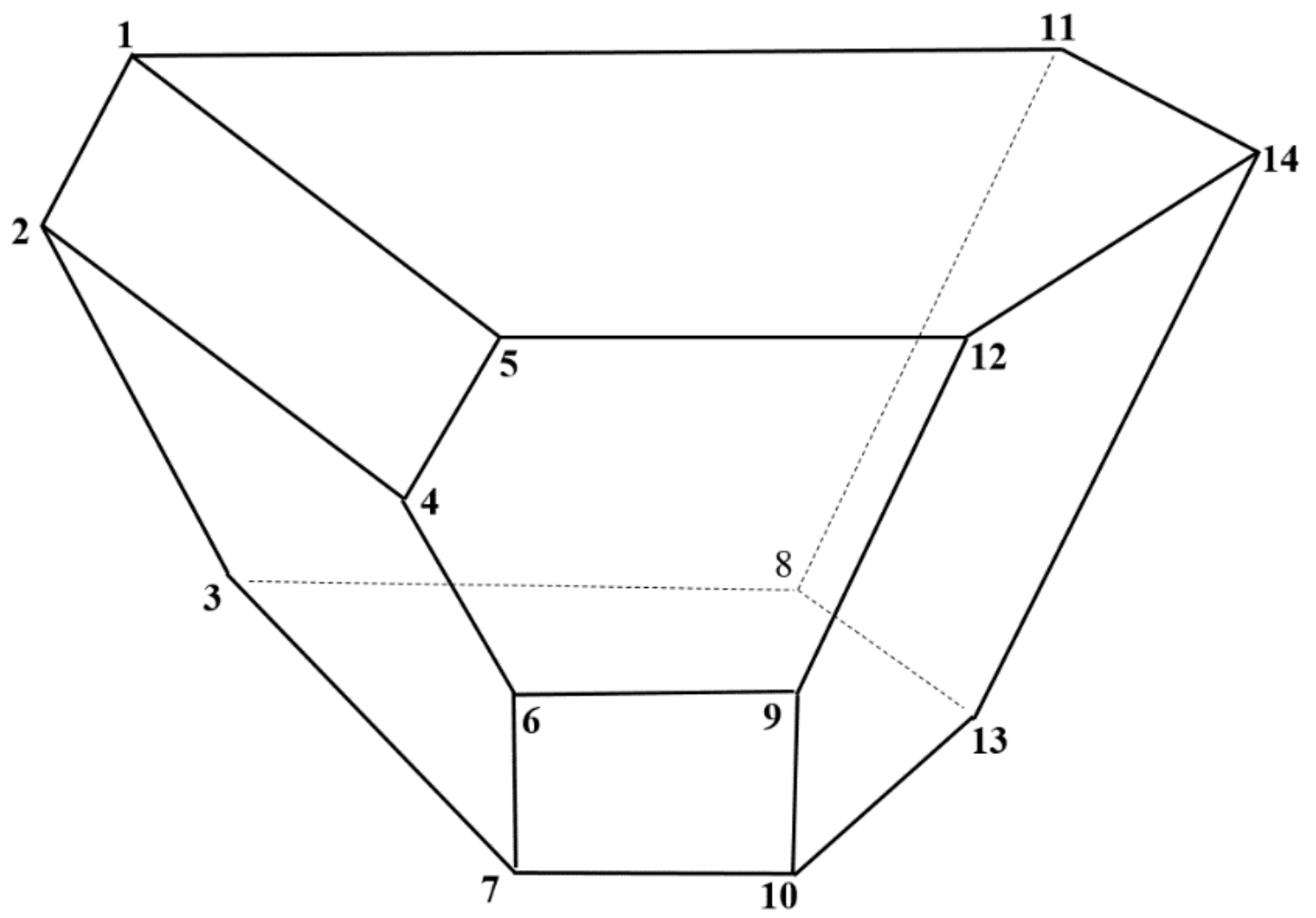

**Figure 1.** The three-dimensional Stasheff polytope

Suppose we are given the facet-defining inequalities shown in Table 2 for describing the feasible space of the above polytope and asked to find unidentified facet-defining inequalities (if any) for this polytope.

**Table 1.** The vertices of the Stasheff polytope.

| Labels | Coordinates |
|---|---|
| 1 | $s_1 = [1\ 6\ 2\ 1]^T$ |
| 2 | $s_2 = [1\ 6\ 1\ 2]^T$ |
| 3 | $s_3 = [1\ 4\ 1\ 4]^T$ |
| 4 | $s_4 = [4\ 3\ 1\ 2]^T$ |
| 5 | $s_5 = [4\ 3\ 2\ 1]^T$ |
| 6 | $s_6 = [4\ 2\ 1\ 3]^T$ |
| 7 | $s_7 = [3\ 2\ 1\ 4]^T$ |
| 8 | $s_8 = [1\ 2\ 3\ 4]^T$ |
| 9 | $s_9 = [4\ 1\ 2\ 3]^T$ |
| 10 | $s_{10} = [3\ 1\ 2\ 4]^T$ |
| 11 | $s_{11} = [1\ 2\ 6\ 1]^T$ |
| 12 | $s_{12} = [4\ 1\ 4\ 1]^T$ |
| 13 | $s_{13} = [2\ 1\ 3\ 4]^T$ |
| 14 | $s_{14} = [2\ 1\ 6\ 1]^T$ |

**Table 2.** Given facets of the Stasheff polytope.

| ID | Given Inequalities | Vertices (Labels) lying in the facets |
|---|---|---|
| 1 | $x_3 \geq 1$ | 2, 3, 4, 6, 7 |
| 2 | $x_2 + x_3 \geq 3$ | 6, 7, 9, 10 |
| 3 | $x_2 \geq 1$ | 9, 10, 12, 13, 14 |
| 4 | $x_1 \leq 4$ | 4, 5, 6, 9, 12 |
| 5 | $x_1 \geq 1$ | 1, 2, 3, 8, 11 |
| 6 | $x_1 + x_2 \geq 3$ | 8, 11, 13, 14 |
| 7 | $x_1 + x_2 \leq 7$ | 1, 2, 4, 5 |

The first step of Algorithm 2 requires us to find all equality constraints (if they are not already identified). Thus, we use the ECA-method (Algorithm 1) to find any new equality constraints as follows:

We have $v_{i-1} := s_i - s_1,\ i = 2, \ldots, 14, V := [v_1|v_2|\ldots|v_{13}] =$

$$\begin{bmatrix} 0 & 0 & 3 & 3 & 3 & 2 & 0 & 3 & 2 & 0 & 3 & 1 & 1 \\ 0 & -2 & -3 & -3 & -4 & -4 & -4 & -5 & -5 & -4 & -5 & -5 & -5 \\ -1 & -1 & -1 & 0 & -1 & -1 & 1 & 0 & 0 & 4 & 2 & 1 & 4 \\ 1 & 3 & 1 & 0 & 2 & 3 & 3 & 2 & 3 & 0 & 0 & 3 & 0 \end{bmatrix}$$

Then, a set of all equality constraints, $\tilde{A}x = \tilde{b}$ can be obtained by the equations in (*) as follows:

$$\tilde{A} = \mathbf{null}(V^T)^T = [1\ 1\ 1\ 1], \qquad \tilde{b} := As_1 = 10.$$

Therefore, the only equality constraint for this polytope is $x_1 + x_2 + x_3 + x_4 = 10$. Subsequently, we need to construct the model (§) and then update it by embedding the following deduplication constraints:

$$\begin{aligned}
&-\theta_1 + \theta_2 + \theta_3 + \theta_4 - \theta_5 + \theta_6 + \theta_7 - \theta_8 - \theta_9 - \theta_{10} - \theta_{11} - \theta_{12} - \theta_{13} - \theta_{14} \leq 4,\\
&-\theta_1 - \theta_2 - \theta_3 - \theta_4 - \theta_5 + \theta_6 + \theta_7 - \theta_8 + \theta_9 + \theta_{10} - \theta_{11} - \theta_{12} - \theta_{13} - \theta_{14} \leq 3,\\
&-\theta_1 - \theta_2 - \theta_3 - \theta_4 - \theta_5 - \theta_6 - \theta_7 - \theta_8 + \theta_9 + \theta_{10} - \theta_{11} + \theta_{12} + \theta_{13} + \theta_{14} \leq 4,\\
&-\theta_1 - \theta_2 - \theta_3 + \theta_4 + \theta_5 + \theta_6 - \theta_7 - \theta_8 + \theta_9 - \theta_{10} - \theta_{11} + \theta_{12} - \theta_{13} - \theta_{14} \leq 4,\\
&+\theta_1 + \theta_2 + \theta_3 - \theta_4 - \theta_5 - \theta_6 - \theta_7 + \theta_8 - \theta_9 - \theta_{10} + \theta_{11} - \theta_{12} - \theta_{13} - \theta_{14} \leq 4,\\
&-\theta_1 - \theta_2 - \theta_3 - \theta_4 - \theta_5 - \theta_6 - \theta_7 + \theta_8 - \theta_9 - \theta_{10} + \theta_{11} - \theta_{12} + \theta_{13} + \theta_{14} \leq 3,\\
&+\theta_1 + \theta_2 - \theta_3 + \theta_4 + \theta_5 - \theta_6 - \theta_7 - \theta_8 - \theta_9 - \theta_{10} - \theta_{11} - \theta_{12} - \theta_{13} - \theta_{14} \leq 3.
\end{aligned}$$

Note that deduplication constraints are obtained by writing the constraint (‡) for each of the existing inequalities given in Table 2. This set of constraints ensures that we do not generate any of existing inequalities.

**Iteration 1.**

Solving (§) by setting $M = 100$ and $\varepsilon = 0.01$, we arrive at the following solution: $\pi = [0.001\ 0.001\ 0.001 - 0.009]^T$, $\pi_0 = 0$, $\theta = [1\ 0\ 0\ 0\ 1\ 0\ 0\ 0\ 0\ 0\ 1\ 1\ 0\ 1]^T$.

Checking the dimension of the convex hull of the extreme points lying on $0.001x_1 + 0.001x_2 + 0.001x_3 - 0.009x_4 \leq 0$ implies that it is a facet-defining inequality. Furthermore, checking (†) verifies that it is associated with a new facet. As mentioned in Remark 1, when a polytope is not full-dimension, a facet can be represented by several inequalities, and one might be interested in finding a simple version of the inequality to represent that facet. For doing so, one can utilize the instruction in Remark 1 and reach $x_4 \leq 1$.

Now we need to add the deduplication constraint to eliminate the currently found facet-defining inequality from the solution space of (§). This constraint is of the form:

$$\theta_1 - \theta_2 - \theta_3 - \theta_4 + \theta_5 - \theta_6 - \theta_7 - \theta_8 - \theta_9 - \theta_{10} + \theta_{11} + \theta_{12} - \theta_{13} + \theta_{14} \leq 4,$$

and, the Constraint (3) of (§) must be updated from $\sum_{i=1}^{N} \theta_i \leq 13$ to $\sum_{i=1}^{N} \theta_i \leq 5$.

**Iteration 2.**

Solving the updated (§), we find the following solution:

$\pi = [-0.004 - 0.004 - 0.004\ \ 0.006]^T$, $\qquad \pi_0 = 0$, $\qquad \theta = [0\ 0\ 1\ 0\ 0\ 0\ 1\ 1\ 0\ 1\ 0\ 0\ 1\ 0]^T$. Dimension and non-redundancy checking reveals that the found inequality is a new facet-defining inequality and by implementing the Remark 1 instruction the simpler version of this inequality can be found in the form of $x_1 + x_2 + x_3 \geq 6$. Hence, this is the second new facet-defining inequality.

The deduplication constraint for this iteration is:

$$-\theta_1 - \theta_2 + \theta_3 - \theta_4 - \theta_5 - \theta_6 + \theta_7 + \theta_8 - \theta_9 + \theta_{10} - \theta_{11} - \theta_{12} + \theta_{13} - \theta_{14} \leq 4,$$

and Constraint (3) of (§) remain intact as this facet-defining constraint is also satisfied by 5 extreme points.

**Iteration 3.**

We terminate Algorithm 2 since the updated (§) model is now infeasible. Therefore, the found facet-defining inequalities are the ones presented in Table 3.

**Table 3.** New identified facets of the Stasheff polytope.

| Vertices (Labels) lying in the facet | Facet-defining Inequalities |
|---|---|
| 1, 5, 11, 12, 14 | $x_4 \le 1$ |
| 3, 7, 8, 10, 13 | $x_1 + x_2 + x_3 \ge 6$ |

Note that facets represented in Table 2 and Table 3 fully describe the Stasheff polytope. Having illustrated the ICA-method with a simple example, we now show how it can be used to strengthen formulations of interesting integer programming problems, such as the Traveling Salesman Problem.

## 4. A New TSP Formulation

In this section, we will first develop a new TSP formulation and then utilize the proposed method in the preceding section to strengthen the formulation by devising some strong inequalities.

The now classical Traveling Salesman Problem (TSP) constitutes a famous challenge for operations researchers, mathematicians, and computers scientists. In particular, suppose there are $n$ cities, and a traveling salesman starts from a home city, passing through all other cities exactly once before returning to the home city. Such a travel path is called a tour or a Hamiltonian cycle (HC). The distance between each pair of cities $i$ and $j$ is given as $C_{ij}$, and so for any tour, the tour length is the sum of distances traveled. Hence, the TSP can be simply thought of as the optimization problem of identifying the tour of shortest length. When $C_{ij} \neq C_{ji}$ the problem is referred to as the Asymmetric Traveling Salesman Problem (ATSP).

We propose a new model for the TSP based on its relationship to the Hamiltonian Cycle Problem (HCP) and corresponding polytope constructions. Feinberg [11] investigated the relationship between the HCP and discounted Markov Decision Processes. He constructed a new polytope corresponding to a given graph $G$, which we shall refer to as $\mathcal{H}_\beta$. The constraints of $\mathcal{H}_\beta$ are of the form:

$$
\begin{aligned}
&\textstyle\sum_{j=2}^{n} x_{1j} - \beta \sum_{j=2}^{n} x_{j1} = 1 - \beta^n \\
&\textstyle\sum_{\substack{j=1 \\ j\neq i}}^{n} x_{ij} - \beta \sum_{\substack{j=1 \\ j\neq i}}^{n} x_{ji} = 0, && i = 2, \dots, n \\
&\textstyle\sum_{j=2}^{n} x_{1j} = 1 && && (\mathcal{H}_\beta) \\
&x_{ij} \ge 0, && i, j = 2, \dots, n, \;\; i \neq j \\
&x_{ij} = 0, && \text{arc}(i, j) \notin G
\end{aligned}
$$

where $\beta$ is a constant value in (0,1). Feinberg showed that, if the graph $G$ is Hamiltonian, then the polytope $\mathcal{H}_\beta$ has an extreme point for each of its Hamiltonian cycles. This result is elaborated in Theorem 5 below.

**Theorem 5.** *If in a feasible solution of $\mathcal{H}_\beta$, for each node $i \in \{1, \dots, n\}$, exactly one variable $x_{ij}$ takes a positive value (i.e., from each node only one flow emanates), then that solution is an extreme point of $\mathcal{H}_\beta$ and the arcs associated with the positive valued variables in that solution trace out a tour or Hamiltonian cycle of the graph $G$.*

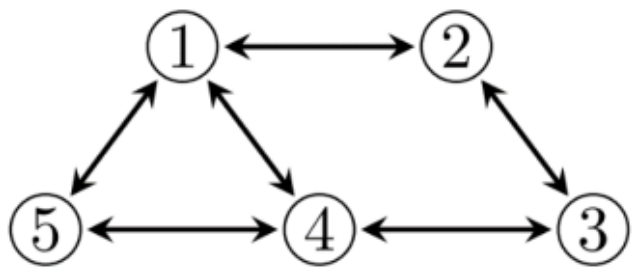

**Figure 2.** A Hamiltonian graph

**Example 2.** The $\mathcal{H}_\beta$ set of constraints for the graph given in Figure 2 can be written as:

$$
\begin{aligned}
&x_{12} + x_{14} + x_{15} - \beta(x_{21} + x_{41} + x_{51}) = 1 - \beta^5 \\
&x_{21} + x_{23} - \beta(x_{12} + x_{32}) = 0 \\
&x_{32} + x_{34} - \beta(x_{23} + x_{43}) = 0 \\
&x_{41} + x_{43} + x_{45} - \beta(x_{14} + x_{34} + x_{54}) = 0 \\
&x_{51} + x_{54} - \beta(x_{15} + x_{45}) = 0 \\
&x_{12} + x_{14} + x_{15} = 1 \\
&x_{12}, x_{14}, x_{15}, x_{21}, x_{23}, x_{32}, x_{34}, x_{41}, x_{43}, x_{45}, x_{51}, x_{54} \geq 0.
\end{aligned}
$$

It can be seen by inspection that

$$
\begin{cases}
x_{12} = 1, x_{23} = \beta, x_{34} = \beta^2, x_{45} = \beta^3, x_{51} = \beta^4 \\
x_{14} = x_{15} = x_{21} = x_{32} = x_{41} = x_{43} = x_{54} = 0
\end{cases}
$$

is an extreme point of $\mathcal{H}_\beta$, and the set of arcs {(1, 2), (2, 3), (3, 4), (4, 5), (5, 1)}, corresponding to the non-zero variables, traces out a Hamiltonian cycle in that graph.

Exploiting the result in Theorem 5, one can utilize binary variables to develop a new TSP formulation. In particular, let $z_{ij} > 0$ if variable $x_{ij} > 0$ and otherwise $z_{ij} = 0$.Then it is easy to write a set of constraints to ensure that only one flow emanates from each node and this guarantees that the feasible space of the model reduces to only Hamiltonian Cycles (tours). Therefore, the new TSP formulation, which we shall refer to as the TSP$_H(\beta)$ formulation, can be written as follows.

$$
\begin{aligned}
\text{Min} \quad & \sum_{i=1}^{n} \sum_{\substack{j=1 \\ j \neq i}}^{n} c_{ij} z_{ij} \\
\text{s.t.} \quad & \\
& \sum_{j=2}^{n} x_{1j} - \beta \sum_{j=2}^{n} x_{j1} = 1 - \beta^n && (1) \\
& \sum_{\substack{j=1 \\ j \neq i}}^{n} x_{ij} - \beta \sum_{\substack{j=1 \\ j \neq i}}^{n} x_{ji} = 0, \quad i = 2, \dots, n && (2) \\
& \sum_{j=2}^{n} x_{1j} = 1 && (3) \\
& x_{1j} = z_{1j}, \quad j = 2, \dots, n && (4) \quad (\text{TSP}_H(\beta)) \\
& x_{j1} = \beta^{n-1} z_{j1}, \quad j = 2, \dots, n && (5)
\end{aligned}
$$

$$\beta^{n-2} z_{ij} \leq x_{ij} \leq \beta z_{ij}, \qquad i,j = 2,\dots,n, \;\; i \neq j \qquad (6)$$

$$\textstyle\sum_{\substack{j=1\\ j\neq i}}^{n} z_{ij} = 1, \qquad i = 1,\dots,n \qquad (7)$$

$$x_{ij} \geq 0, \qquad i,j = 2,\dots,n, \;\; i \neq j \qquad (8)$$

$$z_{ij} \in \{0,1\}, \qquad i,j = 2,\dots,n, \;\; i \neq j \qquad (9)$$

where $z_{ij} = 1$ if arc $(i,j)$ belongs to the tour and otherwise $z_{ij} = 0$. In $\text{TSP}_{\text{H}}(\beta)$ formulation, constraints (1)-(3) are $\mathcal{H}_\beta$ constraints for complete graph and constraints (4)-(7) provide the required condition in theorem 5. Consequently, the only feasible solutions of $\text{TSP}_{\text{H}}(\beta)$ are Hamiltonian cycles and, therefore, $\text{TSP}_{\text{H}}(\beta)$ is a formulation for solving the ATSP. This formulation can be considered as a node-oriented TSP formulation since it contains a set of variables ($x_{ij}$'s) which show the rank order that each node is visited in a tour. For instance, in the tour of example 2, we have $x_{34} = \beta^2$, and node 3 is visited in the third place with node 1 being assigned a rank of one. Generalizing, if $x_{ij} = \beta^{k-1}, k = 1,\dots,n$, it implies that if the tour starts from node 1, we will visit node $i$ in $k^{\text{th}}$ order in that tour.

To the authors' knowledge, there are a few node-oriented TSP models with (at most) quadratic order ($\mathcal{O}(n^2)$) in the number of variables and constraints. The first such model is well-known as the MTZ formulation, which was introduced by Miller, Tucker, and Zemlin in 1960 [23]. The MTZ formulation was improved by Desrochers and Laporte in 1991 [8] through lifting some of its inequality constraints into facets of the underling polytope. In 2002, Sherali and Driscoll [28] proposed a reformulation of the MTZ constraints called the SD model and proved that this model provides an even tighter LP-relaxation than of that the DL model. More recently, Oncan et al. showed that the LP-relaxation of SD is also stronger than that of the well-known single commodity flow model [29]. Furthermore, based on the analysis in [29] and [27], the SD formulation possesses the tightest relaxation among all known TSP formulations with at most ($\mathcal{O}(n^2)$) variables and constraints. This summary of existing literature motivates us to select the SD model for comparison with our newly developed TSP formulation. We shall compare the lower bounds generated by solving the LP-relaxation of both models (SD and $\text{TSP}_{\text{H}}(\beta)$) over a set of ATSP instances from TSPLIB [26]. We will then tighten the $\text{TSP}_{\text{H}}(\beta)$ formulation using the ICA-method developed in earlier sections of this paper to show its effectiveness.

The SD model [28] is defined as:

$$\text{Min} \;\; \textstyle\sum_{i=1}^{n} \sum_{\substack{j=1\\ j\neq i}}^{n} c_{ij} z_{ij}$$

s.t.

$$\textstyle\sum_{j=2}^{n} x_{ij} + (n-1)x_{i1} = u_i, \;\; i = 2,\dots,n \qquad (1)$$

$$\textstyle\sum_{i=2}^{n} x_{ij} + 1 = u_j, \;\; j = 2,\dots,n \qquad (2)$$

$$z_{ij} \leq x_{ij} \leq (n-2)x_{ij}, \;\; i,j = 2,\dots,n, \;\; i \neq j \qquad (3)$$

$$u_j + (n-2)z_{ij} - (n-1)(1 - z_{ji}) \leq x_{ij} + x_{ji} \leq u_j - (1 - z_{ji}), \;\; i,j = 2,\dots,n, \;\; i \neq j \qquad (4)$$

$$2 - x_{1j} + (n-3)x_{j1} \leq u_j \leq (n-2) - (n-3)x_{1j} + x_{j1}, \;\; j = 2,\dots,n \qquad (5) \qquad \text{(SD)}$$

$$\sum_{\substack{j=1 \\ j \neq i}}^{n} z_{ij} = 1, \quad i = 1, \dots, n \tag{6}$$

$$\sum_{\substack{i=1 \\ i \neq j}}^{n} z_{ij} = 1, \quad j = 1, \dots, n \tag{7}$$

$$x_{ij} \geq 0, \quad i, j = 2, \dots, n, \quad i \neq j \tag{8}$$

$$z_{ij} \in \{0,1\}, \quad i, j = 2, \dots, n, \quad i \neq j \tag{9}$$

where $u_j$ denotes the rank order that node $j$ is visited with the initial city being assigned a rank of zero. Also, $z_{ij}$ takes on a value of 1 if the tour transitions from node $i$ to node $j$, and 0 otherwise.

Table 4 presents the result of a computational comparison of the LP-relaxation of SD and $\text{TSP}_\text{H}(\beta)$ formulations. Note that the $\text{TSP}_\text{H}(\beta)$ formulation depends on parameter $\beta$. Interestingly, the structure of the polytope constructed as the LP-relaxation of $\text{TSP}_\text{H}(\beta)$ could vary based $\beta$. Thus, in our numerical study, we consider this parameter very close to one as such values provide greater lower bounds than that of smaller $\beta$'s. This observation is consistent with the results in [9] and [10] in which the extreme points of $\mathcal{H}_\beta$ and another relevant polytope are analytically and computationaly studied.

**Table 4.** LP-relaxation comparison study between $\text{TSP}_\text{H}(\beta)$ and SD formulation over ATSP instances of TSPLIB.

| Problem | $\text{TSP}_\text{H}(0.999)$ | $\text{TSP}_\text{H}(0.9999)$ | SD |
|---|---|---|---|
| br17 | 0.0961 | 0.0997 | 27.6786 |
| ft_53 | 5881.7105 | 5926.0611 | 6118.4042 |
| ft_70 | 37875.7125 | 37968.3275 | 38364.5522 |
| ftv_33 | 1179.1479 | 1184.4144 | 1224.5043 |
| ftv_35 | 1374.6202 | 1380.368 | 1415.5116 |
| ftv_38 | 1431.5212 | 1437.3586 | 1480.0553 |
| ftv_44 | 1511.1038 | 1520.019 | 1573.75 |
| ftv_47 | 1420.5586 | 1433.5657 | 1727.2078 |
| ftv_53 | 1374.6202 | 1380.368 | 1415.5116 |
| ftv_55 | 1420.5586 | 1433.5657 | 1513.2704 |
| ftv_64 | 1696.4783 | 1718.5379 | 1765.3919 |
| ftv_70 | 1746.512 | 1764.0605 | 1859.5775 |
| ftv_170 | 2527.8331 | 2620.6339 | 2698.6786 |
| Kro124p | 33710.3188 | 33951.0923 | 35059.5825 |
| P43 | 147.8366 | 147.9874 | 864.581 |
| Rbg323 | 938.865 | 1281.3476 | 1326 |
| Rbg358 | 781.3627 | 1101.7792 | 1163 |
| Rbg403 | 1375.251 | 2352.9368 | 2465 |
| Rbg443 | 1516.6781 | 2585.8968 | 2720 |
| Ry48p | 12489.6927 | 12514.2894 | 13820.4333 |

Although the numerical results in Table 4 reveal that the lower bound found by $\text{TSP}_\text{H}(0.9999)$ are greater than of that $\text{TSP}_\text{H}(0.999)$, the performance of $\text{TSP}_\text{H}(\beta)$ is far weaker than of that SD in all tested instances. Thus, one might be interested in employing ICA-method to strengthen $\text{TSP}_\text{H}(\beta)$ formulation by generating some double index facet-defining inequalities.

Note that we limited ourselves to find only double index inequalities as we do not wish to exceed using $\mathcal{O}(n^2)$ variables and constraints.

## 5. Identification of Strong Valid Inequalities for TSPH($\boldsymbol{\beta}$)

Before implementing the ICA-method, it is beneficial to conduct symmetry break preparations which can dramatically reduce the required computations. In particular, we set the zero value for many decision variables of the model (§), and this makes Algorithm 2 much more efficient as this model must be solved in each iteration of the algorithm. The symmetry that arises in (§) is due to the indistinguishability of nodes $2, \dots, n$ in the TSPH($\beta$) formulation. Therefore, for any solution of (§), there exist equivalent symmetric reflections of this solution obtained by relabelling the nodes through mapping labels $2, \dots, n$ to any reordering of these nodes (in $(n-1)!$ ways). Note that not necessarily every map result in a new solution. To elaborate what symmetry means here, consider a double index inequality such as constraint (6) in TSPH($\beta$) which can be written for $\binom{n-1}{2}$ different combinations of $i$ and $j$. However, if these inequalities were unknown, we would only need to find one such an inequality to discover that family of inequalities. To prevent finding different members of an inequality family several times, we can permit only one (or a few of them) to be in the feasible space of (§). More precisely, we could assume that the inequality $\pi^T x \leq \pi_0$ that we are trying to find through solving (§) is of the form:

$$\begin{aligned} &\sum_{k=2}^{n} \bar{\pi}_{1k} z_{1k} + \sum_{k=2}^{n} \bar{\pi}_{k1} z_{k1} + \sum_{\substack{k=1\\k\neq 2}}^{n} \bar{\pi}_{2k} z_{2k} + \sum_{\substack{k=1\\k\neq 2}}^{n} \bar{\pi}_{k2} z_{k2} \\ &+ \sum_{\substack{k=1\\k\neq 3}}^{n} \bar{\pi}_{3k} z_{3k} + \sum_{\substack{k=1\\k\neq 3}}^{n} \bar{\pi}_{k3} z_{k3} + \sum_{\substack{k=1\\k\neq 2}}^{n} \pi_{2k} x_{2k} + \sum_{\substack{k=1\\k\neq 2}}^{n} \pi_{k2} x_{k2} \\ &+ \sum_{\substack{k=1\\k\neq 3}}^{n} \pi_{3k} x_{3k} + \sum_{\substack{k=1\\k\neq 3}}^{n} \pi_{k3} x_{k3} \leq \pi_0 \end{aligned} \qquad (**)$$

This is due to distinguishability of node 1 from other nodes in terms of the TSPH($\beta$) formulation viewpoint which allows considering this node as a solo member group. On the other hand, nodes $2, \dots, n$ are indistinguishable in the TSPH($\beta$) formulation and constitute the second group. That means, if we swap the labels of two distinct nodes from the node set $\{2, \dots, n\}$ in the TSPH($\beta$) formulation, then the lower bound generated by the LP-relaxation of TSPH($\beta$) remains intact, however, if we exchange the label of node 1 with any other nodes, the lower bound generated by LP-relaxation could vary. It follows that the most general form of a double index inequality comprises three indices, namely 1, $i$, and $j$ where indices $i$ and $j$ are two distinct arbitrary indices from the node set $\{2, \dots, n\}$. Thus, we just need to allow the variables associated with arcs connected to nodes 1, $i$, and $j$ be involved in $\pi^T x \leq \pi_0$. Without loss of generality, $i$ and $j$ are set to 2 and 3 respectively in (**) and the coefficients of all other variables set to zero.

One might think that the terms $\sum_{k=2}^{n} \pi_{1k} x_{1k} + \sum_{k=2}^{n} \pi_{k1} x_{k1}$ must also be added to the left-hand side of (**). However, it should be noted that the constraints (4) and (5) in TSPH($\beta$) allow us to remove those terms. More precisely, variables $x_{1k}$ and $x_{k1}, k = 2, \dots, n$ can be substituted by $z_{1k}$ and $\beta^{n-1} z_{k1}, k = 2, \dots, n,$ respectively, and as latter variables are involved in $\pi^T x \leq$ $\pi_0$, the former ones can be removed. The elimination of these terms in (**) reduce the symmetries in the solution space of (§) and consequently make it easier to solve. We refer the readers to [19] and [22] for further discussion about symmetry breaking techniques in integer linear programs.

Before any attempt to find strong valid inequalities, we first implemented the ECA-method (Algorithm 1) and found $\sum_{\substack{i=1 \\ i\neq j}}^{n} z_{ij} = 1,\ j = 1, \dots, n$ as the only set of unidentified non-redundant equality constraints. Obviously, this set of constraints along with constraints (7) in TSP$_{\mathrm{H}}(\beta)$ constitute the assignment constraints which appear in most TSP formulations.

In order to strengthen the TSP$_{\mathrm{H}}(\beta)$ formulation by inequalities, we employed Matlab R2017b and CPLEX 12.8 and implemented the ICA-method (Algorithm 2) with parameters $\beta = 0.999$, $\varepsilon = 0.1$, $M = 100$, for $n = 6$ and $n = 7$, to determine the unknown coefficients of inequality (**). All computational experiments are carried out on a HP Z840 workstation with Intel Xeon E5-2640 2.4 GHz 10-Core Processors, 128GB RAM, running Windows 7 operating system and using a single thread.

After completion of this process, we found $\theta$ vectors for facet-defining inequalities. Recall that the $\theta$ vector for an inequality determines which solutions lie on that inequality. Utilizing the instruction in Remark 1, we used $\theta$ vectors and reobtained the inequalities but this time we considered $\beta$ as a parameter and used symbolic Matlab. This leads to finding expressions for the inequalities based on parameter $\beta$ which allows us to explore whether there are any patterns that are generalizable as $n$ grows. By recognizing generalizable patterns in the generated inequalities, we could find a rather large number of valid inequalities, that are facet-defining for those dimensions ($n = 6, 7$). These inequalities, which we shall refer to as Set 3, are listed below:

$$\beta^{n-1}\left(1 - z_{ij} - z_{ji}\right) \le \sum_{\substack{a=1 \\ a\neq i,j}}^{n} x_{ja} - \beta x_{ij} \le \beta\left(1 - z_{ij} - z_{ji}\right), \qquad i,j = 2, \dots, n,\ \ i \neq j \tag{1}$$

$$\beta^{n-2}\left(1 - z_{j1} - z_{1j}\right) + \beta z_{1j} + B^{n-1} z_{j1} \le \sum_{\substack{a=1 \\ a\neq j}}^{n} x_{ja} \le \beta^{2}\left(1 - z_{j1} - z_{1j}\right) + \beta z_{1j} + \beta^{n-1} z_{j1}, \qquad i = 2, \dots, n \tag{2}$$

$$x_{ji} \le (\beta - \beta^{2}) z_{1j} + \beta^{2} z_{ji}, \qquad i,j = 2, \dots, n,\ \ i \neq j \tag{3}$$

$$(\beta^{n-2} - \beta) z_{1j} + \frac{-\beta^{2n-6} + \beta^{n-3}}{\sum_{a=0}^{n-5} \beta^{a}} z_{i1} + \beta^{n-2} z_{j1} - \frac{\beta^{n-3}}{\sum_{a=0}^{n-5} \beta^{a}} z_{ji} + \frac{\sum_{a=0}^{n-4} \beta^{a}}{\sum_{a=0}^{n-5} \beta^{a}} x_{ji} + \sum_{\substack{a=2 \\ a\neq i,j}}^{n} x_{ja} \ge B^{n-2}, \qquad i,j = 2, \dots, n,\ \ i \neq j \tag{4}$$

$$(\beta^{n-3} - \beta) z_{1j} + (\beta^{n-2} - \beta^{n-3}) z_{i1} + \beta^{n-2} z_{j1} + \frac{\beta^{2n-6} + \beta^{n-3}}{\sum_{a=0}^{n-4} \beta^{a}} z_{ji} + \frac{-\beta^{n-4} + \sum_{a=0}^{n-5} \beta^{a}}{\sum_{a=0}^{n-4} \beta^{a}} x_{ji} + \sum_{\substack{a=2 \\ a\neq i,j}}^{n} x_{ja} \ge \beta^{n-2}, \qquad i,j = 2, \dots, n,\ \ i \neq j \tag{5}$$

(Set 3)

$$(\beta^{n-4}-\beta)\left(z_{1i}+z_{1j}\right)+(\beta^{n-2}-\beta^{n-3}+\beta^{n-4})\left(z_{i1}+z_{j1}\right)+ \tag{6}$$
$$\frac{1}{\beta}\left(x_{ij}+x_{ji}\right)+\sum_{\substack{a=2\\a\neq i,j}}^{n} x_{ia}+\sum_{\substack{a=2\\a\neq i,j}}^{n} x_{ja} \geq$$
$$\beta^{n-2}+B^{n-4}, \quad i,j=2,\dots,n,\ i\neq j$$

$$(\beta^{n-4}-\beta)z_{1i}+(\beta^{n-2}-\beta^{n-3}+\beta^{n-4})\left(z_{i1}+z_{j1}\right)+(\beta^{n-4}+\beta^{n-2})z_{ji} \tag{7}$$
$$+\frac{1}{\beta}x_{ij}-\beta x_{ji}+\sum_{\substack{a=2\\a\neq i,j}}^{n} x_{ia}+\sum_{\substack{a=2\\a\neq i,j}}^{n} x_{ja} \geq$$
$$\beta^{n-2}+\beta^{n-4}, \quad i,j=2,\dots,n,\ i\neq j$$

$$(\beta^{n-4}-\beta)\left(z_{1i}+z_{1j}\right)+(\beta^{n-2}-\beta^{n-3}+\beta^{n-4})\left(z_{i1}+z_{j1}\right) \tag{8}$$
$$+\beta^{n-4}\left(\frac{1+\sum_{a=2}^{n-3}\beta^{a}}{\sum_{a=0}^{n-5}\beta^{a}}\right)z_{ji}+\frac{1}{\beta}x_{ij}+\frac{1-\beta^{n-4}-\beta^{n-5}}{\sum_{a=0}^{n-5}\beta^{a}}x_{ji}+\sum_{\substack{a=2\\a\neq i,j}}^{n} x_{ia}$$
$$+\sum_{\substack{a=2\\a\neq i,j}}^{n} x_{ja} \geq \beta^{n-2}+\beta^{n-4}, \quad i,j=2,\dots,n,\ i\neq j$$

$$(\beta^{n-3}-\beta^{2})z_{1i}+B^{n}z_{i1}+(\beta^{n}+\beta^{n-1}+\beta^{n-3})z_{j1}+(\beta^{n}+\beta^{n-1} \tag{9}$$
$$+\beta^{n-2}+\beta^{n-3})z_{ji}+x_{ij}+\beta\sum_{\substack{a=2\\a\neq i,j}}^{n} x_{ia}+(\beta^{2}+\beta+1)\sum_{\substack{a=2\\a\neq i,j}}^{n} x_{ja}$$
$$-\beta^{2}x_{ji} \geq \beta^{n}+\beta^{n-1}+\beta^{n-2}+\beta^{n-3}, \quad i,j=2,\dots,n,\ i\neq j$$

$$\frac{\beta^{n-3}-\beta^{2}}{\beta+1}z_{1i}+(\beta^{n-3}-\beta)z_{1j}+\frac{\beta^{n-1}}{\beta+1}z_{i1}+\frac{\beta^{n-1}+\beta^{n-3}}{\beta+1}z_{j1} \tag{10}$$
$$+\beta^{n-3}\left(\frac{1+\beta+\sum_{a=3}^{n-2}\beta^{a}}{1+\beta^{n-3}+2\sum_{a=1}^{n-4}\beta^{a}}\right)z_{ji}+\frac{1}{B+1}x_{ij}+\frac{1-\beta^{n-3}+2\beta+\sum_{a=2}^{n-5}\beta^{a}}{1+\beta^{n-3}+2\sum_{a=1}^{n-4}\beta^{a}}x_{ji}$$
$$+\frac{\beta}{\beta+1}\sum_{\substack{a=2\\a\neq i,j}}^{n} x_{ia}+\sum_{\substack{a=2\\a\neq i,j}}^{n} x_{ja} \geq \frac{\beta^{n-1}+\beta^{n-2}+\beta^{n-3}}{\beta+1}, \quad i,j=2,\dots,n,\ i\neq j$$

The following set of inequalities hold only for even $n$:

$$\frac{\beta^{n-3}-\beta^{2}}{\beta+1}z_{1i}+(\beta^{n-3}-\beta)z_{1j}+\frac{\beta^{n-1}}{\beta+1}z_{i1}+\frac{\beta^{n-1}+\beta^{n-3}}{\beta+1}z_{j1} \tag{11}$$
$$+\beta^{n-3}\left(\frac{1+\sum_{a=1}^{\frac{n}{2}-2}\beta^{2a+1}}{\sum_{a=0}^{n-4}\beta^{a}}\right)z_{ji}+\frac{1}{\beta+1}x_{ij}+\frac{1-\beta^{n-4}+\sum_{a=0}^{\frac{n}{2}-3}\beta^{2a+1}}{\sum_{a=0}^{n-4}\beta^{a}}x_{ji}$$
$$+\frac{\beta}{\beta+1}\sum_{\substack{a=2\\a\neq i,j}}^{n} x_{ia}+\sum_{\substack{a=2\\a\neq i,j}}^{n} x_{ja} \geq \frac{\beta^{n-1}+\beta^{n-2}+\beta^{n-3}}{\beta+1}, \quad i,j=2,\dots,n,\ i\neq j$$

Checking the above inequalities for larger $n$*'s*, namely $n=8,9,10,11$, we found them not only valid but also facet-defining for those dimensions as well. These results support our following conjecture.

**Conjecture 1.** *There exists a constant* $\beta_0$ *such that for* $n \geq 6$*,and* $\beta \in (\beta_0, 1)$, *the inequalities in Set 3, are facet-defining for the* TSPH$(\beta)$ *formulation.*

Conjecture 1 states that when $\beta$ is large enough, the inequalities in Set 3 are facet-defining for TSPH$(\beta)$ formulation. Note that while it is not difficult to show that these inequalities are valid for any $n \geq 6$, proving that they are facet-defining needs further effort and is beyond the purpose of this paper. For the sake of brevity, we supply explanations for the validity of just three inequalities from Set 3, namely, inequalities (1), (7) and (10) in Appendix 1. The others can be explained by derivations that are conceptually similar.

We are now able to conduct a comparison study between the strengthened $TSP_H(\beta)$ and SD formulations. For strengthening the $TSP_H(\beta)$ formulation, inequalities (1)-(3) of Set 3 are appended to the former set of constraints of $TSP_H(\beta)$, and we shall refer to this model as $TSP_H*(\beta)$. Note that the remainder of the newly developed inequalities in Set 3 were not added to the model as we wish to keep the number of constraints almost equivalent to that of the SD formulation. We conducted the computation for $TSP_H*(\beta)$ (for $\beta$ = 0.999 and 0.9999), and the results are summarized in Table 8. The boldly written values show the instances where $TSP_H*(\beta)$ could not beat the SD formulation.

**Table 8.** LP-relaxation comparison study between $TSP_H*(\beta)$ and SD formulation over ATSP instances of TSPLIB.

| Problem | $TSP_H(0.999)$ | $TSP_H(0.9999)$ | SD |
|---|---|---|---|
| br17 | **27.6555** | **27.6763** | 27.6786 |
| ft_53 | 6118.4166 | 6120.3608 | 6118.4042 |
| ft_70 | **38364.2068** | 38365.1167 | 38364.5522 |
| ftv_33 | 1224.648 | 1224.7385 | 1224.5043 |
| ftv_35 | 1415.5904 | 1415.6172 | 1415.5116 |
| ftv_38 | 1480.1222 | 1480.1248 | 1480.0553 |
| ftv_44 | 1573.75 | 1573.75 | 1573.75 |
| ftv_47 | 1727.2183 | 1727.2556 | 1727.2078 |
| ftv_53 | 1415.5904 | 1415.6172 | 1415.5116 |
| ftv_55 | **1513.2479** | 1513.3305 | 1513.2704 |
| ftv_64 | 1765.4209 | 1765.4596 | 1765.3919 |
| ftv_70 | 1859.6408 | 1859.6454 | 1859.5775 |
| ftv_170 | **2698.6529** | **2698.6773** | 2698.6786 |
| Kro124p | **35056.9562** | 35060.024 | 35059.5825 |
| P43 | **853.3256** | **863.7552** | 864.581 |
| Rbg323 | 1326 | 1326 | 1326 |
| Rbg358 | 1163 | 1163 | 1163 |
| Rbg403 | 2465 | 2465 | 2465 |
| Rbg443 | 2720 | 2720 | 2720 |
| Ry48p | 13820.5307 | 13820.7076 | 13820.4333 |

Naturally, $TSP_H*(0.9999)$ performed better than $TSP_H*(0.999)$ in all instances. Furthermore, $TSP_H*(0.9999)$ outperforms SD in twelve out of twenty test problems, exactly matches SD in five instances and is outperformed slightly by SD in only three instances. Comparing these results to the results in Table 4 also reveals that inequalities (1)-(3) has considerably strengthened the $TSP_H*(\beta)$ formulation.

## 6. Conclusion

In this paper, we propose a generic method to assist modelers in identifying strong linear valid inequalities that can be used to refine feasible spaces of LP-relaxations of mixed integer programmes. The main purpose of this method is to make the devising of strong valid inequalities for IP problems closer to a methodological process rather than an artistic task.

Although the proposed method is still far from a fully automated process, it can be considered as an initial step for this notoriously difficult task. To illustrate the proposed method, a new TSP formulation is introduced and then the ICA-method is employed to devise strong valid inequalities to strengthen the formulation. Comparison between the relaxation of this formulation and a widely used TSP formulation confirms the effectiveness of the devised strong valid inequalities.

Finally, as a naturally arising topic for future research, one can explore whether there are more properties of facet-defining inequalities that can be expressed as new constraints for (§). This can refine the feasible space of (§) from many valid inequalities that are not facets, and therefore the algorithm will be more efficient and more likely to find strong valid inequalities. Another natural line of research is to exploit the ICA-method to extract strong valid inequalities for tightening IP formulations of many other problems.

## References

[1] Assarf, B., Gawrilow, E., Herr, K., Joswig, M., Lorenz, B., Paffenholz, A., & Rehn, T. (2017). Computing convex hulls and counting integer points with polymake. *Mathematical Programming Computation*, *9*(1), 1-38.

[2] Avis, D., Bremner, D., & Seidel, R. (1997). How good are convex hull algorithms?. *Computational Geometry*, *7*(5-6), 265-301.

[3] Avis, D., & Fukuda, K. (1992). A pivoting algorithm for convex hulls and vertex enumeration of arrangements and polyhedra. *Discrete & Computational Geometry*, *8*(3), 295-313.

[4] Avis, D., & Fukuda, K. (1996). Reverse search for enumeration. *Discrete Applied Mathematics*, *65*(1-3), 21-46.

[5] Bektaş, T., & Gouveia, L. (2014). Requiem for the Miller–Tucker–Zemlin subtour elimination constraints?. *European Journal of Operational Research*, *236*(3), 820-832.

[6] Christof, T., Löbel, A., & Stoer, M. (1997). Porta-polyhedron representation transformation algorithm. *Software package, available for download at http://www. zib. de/Optimization/Software/Porta*.

[7] Conforti, M., Cornuéjols, G., & Zambelli, G. (2014). *Integer programming* (Vol. 271). Berlin: Springer.

[8] Desrochers, M., & Laporte, G. (1991). Improvements and extensions to the Miller-Tucker-Zemlin subtour elimination constraints. *Operations Research Letters*, *10*(1), 27-36.

[9] Eshragh, A., & Filar, J. (2011). Hamiltonian cycles, random walks, and discounted occupational measures. *Mathematics of Operations Research*, *36*(2), 258-270.

[10] Eshragh, A., Filar, J. A., Kalinowski, T., & Mohammadian, S. (2018). Hamiltonian cycles and subsets of discounted occupational measures. *arXiv preprint arXiv:1805.04725*.

[11] Feinberg, E. A. (2000). Constrained discounted Markov decision processes and Hamiltonian cycles. *Mathematics of Operations Research*, *25*(1), 130-140.

[12] Filar, J. A., & Moeini, A. (2015). Hamiltonian cycle curves in the space of discounted occupational measures. *Annals of Operations Research*, 1-18.

[13] Fischetti, M., Lodi, A., & Toth, P. (2003). Solving real-world ATSP instances by branch-and-cut. In *Combinatorial Optimization—Eureka, You Shrink!* (pp. 64-77). Springer, Berlin, Heidelberg.

[14] Fukuda, K., Liebling, T. M., & Margot, F. (1997). Analysis of backtrack algorithms for listing all vertices and all faces of a convex polyhedron. *Computational Geometry*, *8*(1), 1-12.

[15] Gawrilow, E., & Joswig, M. (2000). Polymake: a framework for analyzing convex polytopes. In *Polytopes—combinatorics and computation* (pp. 43-73). Birkhäuser, Basel.

[16] Golub, G. H., & Van Loan, C. F. (2012). *Matrix computations* (Vol. 3). JHU press.

[17] Grötschel, M., Lovász, L., & Schrijver, A. (2012). *Geometric algorithms and combinatorial optimization* (Vol. 2). Springer Science & Business Media.

[18] Karp, R. M. (1972). Reducibility among combinatorial problems. In Complexity of computer computations (pp. 85-103). Springer, Boston, MA.

[19] Liberti, L. (2012). Reformulations in mathematical programming: automatic symmetry detection and exploitation. *Mathematical Programming*, *131*(1-2), 273-304.

[20] lusolZ (2016). lusolZ: Nullspace of sparse matrix via LUSOL or SPQR. *available for download at http://stanford.edu/group/SOL/software/lusolZ.*

[21] Bernhard, K., & Vygen, J. (2008). Combinatorial optimization: Theory and algorithms. *Springer, Third Edition, 2005.*

[22] Margot, F. (2010). Symmetry in integer linear programming. In *50 Years of Integer Programming 1958-2008* (pp. 647-686). Springer, Berlin, Heidelberg.

[23] Miller, C. E., Tucker, A. W., & Zemlin, R. A. (1960). Integer programming formulation of traveling salesman problems. *Journal of the ACM (JACM)*, *7*(4), 326-329.

[24] Moeini, A. (2016). *Approximations of the Convex Hull of Hamiltonian Cycles for Cubic Graphs* (Doctoral dissertation, Flinders University).

[25] Moeini, A. (2017). Identification of unidentified equality constraints for integer programming problems. *European Journal of Operational Research*, *260*(2), 460-467.

[26] Reinelt, G. (1991). TSPLIB—A traveling salesman problem library. *ORSA journal on computing*, *3*(4), 376-384.

[27] Roberti, R., & Toth, P. (2012). Models and algorithms for the asymmetric traveling salesman problem: an experimental comparison. *EURO Journal on Transportation and Logistics*, *1*(1-2), 113-133.

[28] Sherali, H. D., & Driscoll, P. J. (2002). On tightening the relaxations of Miller-Tucker-Zemlin formulations for asymmetric traveling salesman problems. *Operations Research*, *50*(4), 656-669.

[29] Öncan, T., Altınel, İ. K., & Laporte, G. (2009). A comparative analysis of several asymmetric traveling salesman problem formulations. *Computers & Operations Research*, *36*(3), 637-654.

[30] Weyl, H. (1934). Elementare theorie der konvexen polyeder. *Commentarii Mathematici Helvetici*, *7*(1), 290-306.; Kuhn, H. W., & Tucker, A. W. (Eds.). (1950). *Contributions to the Theory of Games* (Vol. 1). Princeton University Press.

[31] Wolsey, L. A. (1989). Strong formulations for mixed integer programming: a survey. *Mathematical Programming*, *45*(1), 173-191.

[32] Wolsey, L. A., & Nemhauser, G. L. (2014). *Integer and combinatorial optimization*. John Wiley & Sons.

[33] Ziegler, G. M. (2012). *Lectures on polytopes* (Vol. 152). Springer Science & Business Media.

## Appendix 1

A generic and easy approach to show the validity of inequalities in Set 3 is to consider all the possibilities that binary variables on those inequalities can create. That is, if an inequality has $k$ binary variables, we need to prove that the inequality is valid for all $2^k$ possibilities. For example, a proof for inequality (1) of Set 3 is as follows:

**Proof.** This constraint has two binary variables which leads to four possibilities. These four possible cases are analyzed in Table 5.

**Table 5.** Analysis of proving Inequality (1) from Set 3.

| Case | $z_{ij}$ | $z_{ji}$ | $\pi^T x \vert z_{ij}, z_{ji}$ | $\min\{\pi^T x \vert z_{ij}, z_{ji}\}$ | $\max\{\pi^T x \vert z_{ij}, z_{ji}\}$ |
|---|---|---|---|---|---|
| 1 | 1 | 1 | NA | NA | NA |
| 2 | 1 | 0 | $\sum_{\substack{a=1 \\ a\neq i,j}}^{n} x_{ja} - \beta x_{ij}$ | 0 | 0 |
| 3 | 0 | 1 | $\sum_{\substack{a=1 \\ a\neq i,j}}^{n} x_{ja} - \beta x_{ij}$ | 0 | 0 |
| 4 | 0 | 0 | $\sum_{\substack{a=1 \\ a\neq i,j}}^{n} x_{ja} - \beta x_{ij}$ | $\beta^{n-1}$ | $\beta$ |

In case 1, we have $z_{ij} = z_{ji} = 1$, which is not applicable (NA) as TSP$_{\text{H}}(\beta)$ is a TSP formulation and therefore its constraints never allow arcs $(i,j)$ and $(j,i)$ to appear simultaneously in a tour. In Case 2, we have $z_{ij} = 1$ and $z_{ji} = 0$. Without loss of generality, suppose node $i$ is the $k^{\text{th}}$ position in a tour, so $x_{ij} = \beta^{k-1}$ and according to the constraint (2) in TSP$_{\text{H}}(\beta)$ we have $\sum_{\substack{a=1 \\ a\neq i,j}}^{n} x_{ja} = \beta^k$ which results in $\sum_{\substack{a=1 \\ a\neq i,j}}^{n} x_{ja} - \beta x_{ij} = 0$. In case 3, $z_{ij} = 0$ and constraints (4)-(6) of TSP$_{\text{H}}(\beta)$ imply $x_{ij} = 0$. Also, $z_{ji} = 1$ which implies $x_{ji} > 0$ and $x_{ja} = 0, a = 1, \dots, n, a \neq i,j$ and therefore $\sum_{\substack{a=1 \\ a\neq i,j}}^{n} x_{ja} = 0$. This leads to $\sum_{\substack{a=1 \\ a\neq i,j}}^{n} x_{ja} - \beta x_{ij} = 0$. In case 4, $x_{ij} = 0$ as $z_{ij} = 0$ and obviously $\beta^{n-1} \leq \sum_{\substack{a=1 \\ a\neq i,j}}^{n} x_{ja} \leq \beta$ so $\beta^{n-1} \leq \sum_{\substack{a=1 \\ a\neq i,j}}^{n} x_{ja} - \beta x_{ij} \leq \beta$ and this completes the proof.

■

Similarly, constraint (7) of Set 3 can be proved as follows:

**Proof.** Note that this inequality is of the form $\pi^T x \geq \pi_0$. Hence, one way to prove the validity of this inequality is to show $\min\{\pi^T x\} \geq \pi_0$. For doing so, we evaluated the inequality (7) based on the values its binary variables could take $(\pi^T x | z_{1i}, z_{1j}, z_{i1}, z_{j1})$ in column 5 of Table 6, and in column 6 of this table, we set the remaining decision variables so that the minimum for $\pi^T x$ is achieved $(\min\{\pi^T x\} | z_{1i}, z_{1j}, z_{i1}, z_{j1})$.

**Table 6.** Analysis of proving Inequality (6) from Set 3.

| Case | $z_{1i}$ | $z_{1j}$ | $z_{i1}$ | $z_{j1}$ | $\pi^T x \vert z_{1i}, z_{1j}, z_{i1}, z_{j1}$ | $\min\{\pi^T x \vert z_{1i}, z_{1j}, z_{i1}, z_{j1}\}$ |
|---|---|---|---|---|---|---|
| 1 | 1 | 1 | 1 | 1 | NA | NA |
| 2 | 1 | 1 | 1 | 0 | NA | NA |
| 3 | 1 | 1 | 0 | 1 | NA | NA |
| 4 | 1 | 1 | 0 | 0 | NA | NA |
| 5 | 1 | 0 | 1 | 1 | NA | NA |
| 6 | 1 | 0 | 1 | 0 | NA | NA |
| 7 | 1 | 0 | 0 | 1 | $\beta^{n-4} + \beta^{n-2} - \beta^{n-3} + \beta^{n-4} + \sum_{\substack{a=2 \\ a\neq i,j}}^{n} x_{ia}$ | $2\beta^{n-4} + \beta^{n-2} - \beta^{n-3} + \beta$ |
| 8 | 1 | 0 | 0 | 0 | $\beta^{n-4} - \beta + \frac{1}{\beta}x_{ij} + \sum_{\substack{a=2 \\ a\neq i,j}}^{n} x_{ia} + \sum_{\substack{a=2 \\ a\neq i,j}}^{n} x_{ja}$ | $\beta^{n-4} + \beta^{n-2}$ |
| 9 | 0 | 1 | 1 | 1 | NA | NA |
| 10 | 0 | 1 | 1 | 0 | $\beta^{n-4} + \beta^{n-2} - \beta^{n-3} + \beta^{n-4} + \sum_{\substack{a=2 \\ a\neq i,j}}^{n} x_{ia}$ | $2\beta^{n-4} + \beta^{n-2} - \beta^{n-3} + \beta$ |
| 11 | 0 | 1 | 0 | 1 | NA | NA |
| 12 | 0 | 1 | 0 | 0 | $\beta^{n-4} - \beta + \frac{1}{\beta}x_{ij} + \sum_{\substack{a=2 \\ a\neq i,j}}^{n} x_{ia} + \sum_{\substack{a=2 \\ a\neq i,j}}^{n} x_{ja}$ | $\beta^{n-4} + \beta^{n-2}$ |
| 13 | 0 | 0 | 1 | 1 | NA | NA |
| 14 | 0 | 0 | 1 | 0 | $\beta^{n-2} - \beta^{n-3} + \beta^{n-4} + \frac{1}{\beta}\left(x_{ij} + x_{ji}\right) + \sum_{\substack{a=2 \\ a\neq i,j}}^{n} x_{ia} + \sum_{\substack{a=2 \\ a\neq i,j}}^{n} x_{ja}$ | $\beta^{n-2} + \beta^{n-4} + \beta^{n-1}$ |
| 15 | 0 | 0 | 0 | 1 | $\beta^{n-2} - \beta^{n-3} + \beta^{n-4} + \frac{1}{\beta}\left(x_{ij} + x_{ji}\right) + \sum_{\substack{a=2 \\ a\neq i,j}}^{n} x_{ia} + \sum_{\substack{a=2 \\ a\neq i,j}}^{n} x_{ja}$ | $\beta^{n-2} + \beta^{n-4} + \beta^{n-1}$ |
| 16 | 0 | 0 | 0 | 0 | $\frac{1}{\beta}\left(x_{ij} + x_{ji}\right) + \sum_{\substack{a=2 \\ a\neq i,j}}^{n} x_{ia} + \sum_{\substack{a=2 \\ a\neq i,j}}^{n} x_{ja}$ | $\beta^{n-4} + \beta^{n-2}$ |

It can be seen by inspection that all values in column 6 are greater than (or equal to) $\beta^{n-2} + \beta^{n-4}$ for any $\beta \in (0,1)$ which shows that this inequality is valid in general form.

■

As a more complicated inequality, we prove the validity of inequality (10) of Set 3, below:

**Proof.** For this inequality we have:

**Table 7.** Analysis of proving Inequality (10) from Set 3.

| Case | $z_{1i}$ | $z_{1j}$ | $z_{i1}$ | $z_{j1}$ | $z_{ji}$ | $\pi^T x \mid z_{1j}, z_{1j}, z_{i1}, z_{j1}, z_{ji}$ | $\min\{\pi^T x \mid z_{1j}, z_{1j}, z_{i1}, z_{j1}, z_{ji}\}$ |
|---|---|---|---|---|---|---|---|
| 1 | 1 | 1 | 1 | 1 | 1 | NA | NA |
| 2 | 1 | 1 | 1 | 1 | 0 | NA | NA |
| 3 | 1 | 1 | 1 | 0 | 1 | NA | NA |
| 4 | 1 | 1 | 1 | 0 | 0 | NA | NA |
| 5 | 1 | 1 | 0 | 1 | 1 | NA | NA |
| 6 | 1 | 1 | 0 | 1 | 0 | NA | NA |
| 7 | 1 | 1 | 0 | 0 | 1 | NA | NA |
| 8 | 1 | 1 | 0 | 0 | 0 | NA | NA |
| 9 | 1 | 0 | 1 | 1 | 1 | NA | NA |
| 10 | 1 | 0 | 1 | 1 | 0 | NA | NA |
| 11 | 1 | 0 | 1 | 0 | 1 | NA | NA |
| 12 | 1 | 0 | 1 | 0 | 0 | NA | NA |
| 13 | 1 | 0 | 0 | 1 | 1 | NA | NA |
| 14 | 1 | 0 | 0 | 1 | 0 | $\frac{\beta^{n-3}-\beta^2}{\beta+1}+\frac{\beta^{n-1}+\beta^{n-3}}{\beta+1}+\frac{\beta}{\beta+1}\sum_{\substack{a=2 \\ a\neq i,j}}^{n} x_{ia}$ | $\frac{2\beta^{n-3}+\beta^{n-1}}{\beta+1}$ |
| 15 | 1 | 0 | 0 | 0 | 1 | NA | NA |
| 16 | 1 | 0 | 0 | 0 | 0 | $\frac{\beta^{n-3}-\beta^2}{\beta+1}+\frac{1}{\beta+1}x_{ij}+\sum_{\substack{a=2 \\ a\neq i,j}}^{n} x_{ja}$ | $\frac{\beta^{n-1}+\beta^{n-2}+\beta^{n-3}}{\beta+1}$ |
| 17 | 0 | 1 | 1 | 1 | 1 | NA | NA |
| 18 | 0 | 1 | 1 | 1 | 0 | NA | NA |
| 19 | 0 | 1 | 1 | 0 | 1 | NA | NA |
| 20 | 0 | 1 | 1 | 0 | 0 | $\beta^{n-3}-\beta+\frac{\beta^{n-1}}{\beta+1}+\sum_{\substack{a=2 \\ a\neq i,j}}^{n} x_{ja}$ | $\frac{\beta^{n-1}+\beta^{n-2}+\beta^{n-3}}{\beta+1}$ |
| 21 | 0 | 1 | 0 | 1 | 1 | NA | NA |
| 22 | 0 | 1 | 0 | 1 | 0 | NA | NA |
| 23 | 0 | 1 | 0 | 0 | 1 | $\beta^{n-3}-\beta+\beta^{n-3}\left(\frac{1+\beta+\sum_{a=3}^{n-2}\beta^a}{1+\beta^{n-3}+2\sum_{a=1}^{n-4}\beta^a}\right)+\frac{1-\beta^{n-3}+2B+\sum_{a=2}^{n-5}\beta^a}{1+\beta^{n-3}+2\sum_{a=1}^{n-4}\beta^a}x_{ji}+\frac{\beta}{\beta+1}\sum_{\substack{a=2 \\ a\neq i,j}}^{n} x_{ia}$ | $\beta^{n-3}-\beta+\beta^{n-3}\left(\frac{1+\beta+\sum_{a=3}^{n-2}\beta^a}{1+\beta^{n-3}+2\sum_{a=1}^{n-4}\beta^a}\right)+\frac{1-\beta^{n-3}+2B+\sum_{a=2}^{n-5}\beta^a}{1+\beta^{n-3}+2\sum_{a=1}^{n-4}\beta^a}\beta+\frac{\beta^3}{\beta+1}$ |
| 24 | 0 | 1 | 0 | 0 | 0 | $\beta^{n-3}-\beta+\frac{\beta}{\beta+1}\sum_{\substack{a=2 \\ a\neq i,j}}^{n} x_{ia}+\sum_{\substack{a=2 \\ a\neq i,j}}^{n} x_{ja}$ | $\frac{\beta^{n-1}+\beta^{n-2}+\beta^{n-3}}{\beta+1}$ |
| 25 | 0 | 0 | 1 | 1 | 1 | NA | NA |
| 26 | 0 | 0 | 1 | 1 | 0 | NA | NA |
| 27 | 0 | 0 | 1 | 0 | 1 | $\frac{\beta^{n-1}}{B+1} + \beta^{n-3}\left(\frac{1+B+\sum_{a=3}^{n-2}\beta^a}{1+\beta^{n-3}+2\sum_{a=1}^{n-4}\beta^a}\right) + \frac{1-\beta^{n-3}+2\beta+\sum_{a=2}^{n-5}\beta^a}{1+\beta^{n-3}+2\sum_{a=1}^{n-4}\beta^a} x_{ji}$ | $\frac{\beta^{n-1}}{\beta+1} + \beta^{n-3}\left(\frac{1+\beta+\sum_{a=3}^{n-2}\beta^a}{1+\beta^{n-3}+2\sum_{a=1}^{n-4}\beta^a}\right) + \frac{1-\beta^{n-3}+2\beta+\sum_{a=2}^{n-5}\beta^a}{1+\beta^{n-3}+2\sum_{a=1}^{n-4}\beta^a}\beta^{n-2}$ |
| 28 | 0 | 0 | 1 | 0 | 0 | $\frac{\beta^{n-1}}{\beta+1} + \sum_{\substack{a=2 \\ a\neq i,j}}^{n} x_{ja}$ | $\frac{\beta^{n-1}+\beta^{n-2}+\beta^{n-3}}{\beta+1}$ |
| 29 | 0 | 0 | 0 | 1 | 1 | NA | NA |
| 30 | 0 | 0 | 0 | 1 | 0 | $\frac{\beta^{n-1}+\beta^{n-3}}{\beta+1} + \frac{1}{\beta+1}x_{ij} + \frac{\beta}{\beta+1}\sum_{\substack{a=2 \\ a\neq i,j}}^{n} x_{ia}$ | $\frac{\beta^{n-1}+\beta^{n-2}+\beta^{n-3}}{\beta+1}$ |
| 31 | 0 | 0 | 0 | 0 | 1 | $\beta^{n-3}\left(\frac{1+\beta+\sum_{a=3}^{n-2}\beta^a}{1+\beta^{n-3}+2\sum_{a=1}^{n-4}\beta^a}\right) + \frac{1-\beta^{n-3}+2\beta+\sum_{a=2}^{n-5}\beta^a}{1+\beta^{n-3}+2\sum_{a=1}^{n-4}\beta^a} x_{ji} + \frac{\beta}{\beta+1}\sum_{\substack{a=2 \\ a\neq i,j}}^{n} x_{ia}$ | $\beta^{n-3}\left(\frac{1+\beta+\sum_{a=3}^{n-2}\beta^a}{1+B^{n-3}+2\sum_{a=1}^{n-4}B^a}\right) + \frac{1-\beta^{n-3}+2\beta+\sum_{a=2}^{n-5}v^a}{1+\beta^{n-3}+2\sum_{a=1}^{n-4}\beta^a}\beta^{n-3} + \frac{\beta^{n-1}}{\beta+1}$ |
| 32 | 0 | 0 | 0 | 0 | 0 | $\frac{1}{\beta+1}x_{ij} + \frac{\beta}{\beta+1}\sum_{\substack{a=2 \\ a\neq i,j}}^{n} x_{ia} + \sum_{\substack{a=2 \\ a\neq i,j}}^{n} x_{ja}$ | $\frac{\beta^{n-1}+\beta^{n-2}+\beta^{n-3}}{\beta+1}$ |

As is shown in the last column of the above table, it is easy to see $\min\{\pi^T x\}\,|z_{1j}, z_{1j}, z_{i1}, z_{j1}, z_{ji} \geq \frac{\beta^{n-1}+\beta^{n-2}+\beta^{n-3}}{\beta+1}$ for all cases except cases 23, 27 and 31. So we further simplified the inequality for these three cases and obtained $\beta^{n-1} \geq \beta^{2n-4}$, $0 \geq 0$, and $\beta^{n-2}+1+\beta \geq \beta^{n-3}+\beta^{n-4}+\beta^2$ respectively. All these three inequalities are obviously correct for any $\beta \in (0,1)$ and therefore inequality (10) is valid in general form.

■